\documentclass[11pt]{amsart}

\usepackage{epsfig, color}

\usepackage{amsmath}
\usepackage{amsthm}
\usepackage{hyperref}
\usepackage[margin=1.3in]{geometry}

\numberwithin{equation}{section}

\newtheorem{prop}{Proposition}
\newtheorem{lemma}[prop]{Lemma}

\newtheorem{thm}[prop]{Theorem}
\newtheorem{cor}[prop]{Corollary}

\numberwithin{prop}{section}

\theoremstyle{definition}
\newtheorem{defn}[prop]{Definition}

\newtheorem{rmk}[prop]{Remark}

\newcommand{\del}{\partial}
\newcommand{\delb}{\bar{\partial}}\newcommand{\dt}{\frac{\partial}{\partial t}}
\newcommand{\brs}[1]{\left| #1 \right|}

\newcommand{\gG}{\Gamma}
\renewcommand{\gg}{\gamma}
\newcommand{\gD}{\Delta}
\newcommand{\gd}{\delta}

\newcommand{\gl}{\lambda}

\newcommand{\gw}{\omega}
\newcommand{\ga}{\alpha}
\newcommand{\gb}{\beta}
\renewcommand{\ge}{\epsilon}
\newcommand{\N}{\nabla}
\newcommand{\FF}{\mathcal F}

\newcommand{\til}[1]{\widetilde{#1}}
\newcommand{\nm}[2]{|| #1 ||_{#2}}

\renewcommand{\bar}[1]{\overline{#1}}

\DeclareMathOperator{\Rc}{Rc}
\DeclareMathOperator{\Rm}{Rm}
\DeclareMathOperator{\inj}{inj}

\DeclareMathOperator{\Id}{Id}

\DeclareMathOperator{\grad}{grad}
\DeclareMathOperator{\Vol}{Vol}

\DeclareMathOperator{\supp}{supp}

\begin{document}

\title[The long time behavior of fourth order curvature flows]{The long time
behavior of fourth order curvature flows}
\address{Rowland Hall\\
         University of California, Irvine\\
         Irvine, CA 92617}
\email{\href{mailto:jstreets@uci.edu}{jstreets@uci.edu}}

\thanks{The author was partly supported by a grant from the National Science
Foundation}

\begin{abstract} We show precompactness results for solutions to parabolic
fourth order
geometric evolution equations.  As part of the proof we obtain smoothing
estimates for these flows in the presence of a curvature bound, an improvement
on prior results which also require a Sobolev constant bound.  As consequences
of these results we show
that for any solution with a finite time singularity, the $L^{\infty}$ norm of
the
curvature must go to infinity.  Furthermore, we characterize the
behavior at infinity of solutions with bounded curvature.
\end{abstract}

\date{\today}

\maketitle

\section{Introduction}

We say that a one-parameter family of metrics $(M^n, g(t))$ is a solution to
\emph{fourth-order curvature flow} (FOCF) if $g(t)$ and the associated
Riemannian curvatures $\Rm(t)$ satisfy
\begin{gather} \label{generalflow}
\begin{split}
\dt g =&\ \N^2 \Rm * g + \Rm^{*2},\\
\dt \Rm =&\ - \gD^2 \Rm + \N^2 \Rm * \Rm + \N \Rm^{*2} + \Rm^{*3}.
\end{split}
\end{gather}
\noindent Here the notation $A * B$ refers to some metric contraction of the
tensor product $A \otimes B$.  Solutions to this system arise naturally as
gradient flows of
quadratic
curvature functionals on Riemannian manifolds.   Before discussing the results
of this paper let us give some natural instances of FOCF.  In particular, with
the convention that $\brs{\Rm}^2 = g^{ip} g^{jq} g^{kr} g^{ls} R_{ijkl}
R_{pqrs}$,
let
\begin{align*}
\mathcal F(g) := \frac{1}{2} \int_M \brs{\Rm}^2_g dV_g.
\end{align*}
A basic calculation
(\cite{Besse} Proposition 4.70) shows that
\begin{align} \label{gradF}
\grad \FF =&\ \gd d \Rc - \check{R} + \frac{1}{4} \brs{\Rm}^2 g.
\end{align}
where $d$ is the exterior derivative acting on the Ricci tensor treated
as a one-form with values in the tangent bundle, and $\gd$ is the adjoint of
$d$.  Moreover, 
\begin{align*}
\check{R}_{ij} = R_{i p q r} R_j^{\ p q r}.
\end{align*}
Suppose $(M^n, g(t))$ is a solution to the negative gradient flow for $\mathcal
F$, i.e.
\begin{gather} \label{L2flow}
\begin{split}
\dt g =&\ - \grad \FF.
\end{split}
\end{gather}
It follows from \cite{StreetsL21} that $g(t)$ is a solution to FOCF.

Next, fix a K\"ahler manifold $(M^{2n}, J, g)$.  One may consider the $L^2$ norm
of the scalar curvature restricted to the given K\"ahler class $[\omega]$. 
Expressing a general metric in the K\"ahler class as $\omega_{\phi} := \omega +
\sqrt{-1} \del \delb \phi$, we have
\begin{align*}
\mathcal C(\phi) := \int_M s_{\phi}^2 dV_{\phi}.
\end{align*}
The gradient of this functional is
\begin{align*}
\grad \mathcal C = \sqrt{-1} \del \delb s.
\end{align*}
The associated gradient flow reduces to an equation for the K\"ahler potential,
called the Calabi flow:
\begin{align*}
\dt \phi =&\ s - \frac{\int_M s dV}{\int_M dV}.
\end{align*}
As is shown in \cite{ChenHe2} the associated one-parameter family of metrics
$g_{\phi_t}$ is a solution to FOCF.

Recently Bour \cite{Bour} studied the gradient flow of more general
quadratic curvature functionals in dimension $4$, obtaining a characterization
of the sphere which improves the main result of \cite{StreetsL23}.  These
equations are not quite of the form (\ref{generalflow}), as the leading order
term in the curvature evolution is slightly different, though still parabolic. 
However, since the key input to the strategy of the proof is the fact that one
has local smoothing estimates in Sobolev spaces, the results in this paper
should carry over to these flows mutatis mutandis.

Finally we
point out that a family of higher order geometric evolutions was recently
introduced in \cite{Bachflow} with an aim toward prescribing the ambient
obstruction tensor of Fefferman-Graham.  In dimension four this flow corresponds
to a certain conformally modified gradient flow for the $L^2$ norm of the Weyl
curvature tensor, and thus is
perhaps appropriately called \emph{Bach flow}.  This flow is also of the form
(\ref{generalflow}).

At this point many analytic results have been established for these examples of
FOCF.  Below we will recall some results pertaining to solutions to
(\ref{L2flow}), as they are indicative generally of the types of results which
are available for these flows.  In \cite{StreetsL21} we established the general
short-time existence
and a long-time existence obstruction for solutions to (\ref{L2flow}). 
Specifically we have
\begin{thm} (\cite{StreetsL21} Theorem 6.2) \label{oldLTE} Let $(M^n, g_0)$ be a
compact Riemannian manifold.  The
solution to (\ref{L2flow}) with initial condition $g_0$ exists on a maximal time
interval $[0, T)$.  Furthermore, if $T < \infty$ then either
\begin{align*}
\limsup_{t \to T} \brs{\Rm}_{g(t)} = \infty
\end{align*}
or
\begin{align*}
\limsup_{t \to T} C_S(g(t)) = \infty
\end{align*}
where $C_S(g(t))$ denotes the $L^2$ Sobolev constant of the time-dependent
metrics.
\end{thm}
\noindent The proof exploits $L^2$ derivative estimates satisfied by solutions
of (\ref{L2flow}) in the
presence of a curvature bound, then uses the Sobolev inequality to produce
pointwise bounds from these estimates.  Long-time existence criteria for the
Calabi flow can be found in \cite{ChenHe1}, \cite{ChenHe2}.  Next we recall a
compactness theorem for solutions to (\ref{L2flow}).

\begin{thm} \label{strongcompactness} (\cite{StreetsL21} Theorem 7.1) Let
$(\{M_i, g_i(t), p_i \})$ be a sequence of pointed solutions to
(\ref{L2flow}) on compact manifolds $M_i$, where $t \in (\ga, \gw), - \infty
\leq
\ga < \gw \leq \infty$.  Suppose there exists $C, K < \infty$ such that
\begin{align*}
\sup_{M_i \times (\ga, \gw)} \brs{\Rm(g_i)}_{g_i} \leq K,
\end{align*}
\begin{align*}
\lim_{t \searrow \ga} \nm{\Rm}{L^2(g_i(t))} \leq C, 
\end{align*}
and $\gd > 0$ such that
\begin{align*}
\inf_{M_i \times (\ga, \gw)} \inj_{g_i} \geq \gd.
\end{align*}
Then there exists a subsequence which converges to a complete pointed
solution to (\ref{L2flow}) $(M_{\infty}, g_{\infty}(t), p_{\infty})$
with the same
bounds on curvature and injectivity radius.
\end{thm}
\noindent A similar result for the Calabi flow is implicit in \cite{ChenHe1}. 
This theorem is analogous to Hamilton's compactness result for Ricci flow
(\cite{Hamilton} Theorem 1.2).  A weakness of using this result to understand
blowup limits of solutions to (\ref{L2flow}) is of course the injectivity radius
assumption.  Observe that in fact we are assuming a uniform lower bound at all
points of the manifold.  This is because this implies a uniform upper bound on
the Sobolev constants of the time evolving manifolds, which the proof requires
since it relies on $L^2$ bounds to obtaining the higher order estimates.  One
does not have a noncollapsing estimate at the scale of
curvature for solutions to (\ref{L2flow}) akin to Perelman's estimate for
solutions to Ricci flow (\cite{Perelman}) therefore at present Theorem
\ref{strongcompactness} cannot be applied to understand general blowup limits of
solutions to FOCF.  

Our first new theorem is a generalization of Theorem
\ref{strongcompactness} to sequences of solutions of FOCF with uniform $C^k$
bounds on curvature but no injectivity radius estimate.  This theorem appears in
section 3, and both the statement and proof are analogous to a weak compactness
theorem proved for solutions to Ricci flow by Glickenstein \cite{Glickenstein}. 
In fact Glickenstein's theorem does not use the precise structure of Ricci flow,
except to exploit the derivative estimates which automatically hold in the
presence of a curvature bound.  For our statement we have assumed the necessary bounds and so
his proof is easy to adapt to this setting.  For solutions to FOCF, direct
smoothing techniques only yield bounds for
curvature in $H_k^2$, and require a Sobolev constant estimate to convert to
pointwise bounds, something we are specifically trying to avoid.  However, in an
interesting twist, we
are able to exploit this compactness theorem to \emph{prove} $C^k$
smoothing estimates without a Sobolev constant bound.  Specifically, we have:

\begin{thm} \label{globalsmoothing} Fix $m, n \geq 0$.  There exists a constant
$C = C(m, n)$ so that if $(M^n, g(t))$ is a complete solution to FOCF on
$\left[0, T \right]$ satisfying
\begin{align*}
\sup_{M \times \left[0, T \right]} \brs{\Rm} \leq K,
\end{align*}
then for all $t \in (0, T]$,
\begin{align} \label{globalsmoothingestimate}
\sup_{M} \brs{\N^m \Rm}_{g(t)} \leq C \left( K +
\frac{1}{t^{\frac{1}{2}}} \right)^{1 + \frac{m}{2}}.
\end{align}
\end{thm}

\begin{rmk} In a recent paper on the Calabi flow on toric varieties,
\cite{Huang}, a bound for the first derivative of curvature is required in the
presence of a curvature bound.  Theorem \ref{globalsmoothing} provides this
estimate.  From there the obstruction to long time existence is given by a
certain bound on the derivative of scalar curvature.
\end{rmk}

We can now turn around and exploit
Theorem \ref{globalsmoothing} and the weak compactness theorem of section 3 to
obtain a stronger compactness result as a corollary.

\begin{cor} \label{wkcompactness}
Let $\{(M^n_i, g_i(t), p_i)\}$ be a sequence of complete pointed
solutions of FOCF, where $t \in (\ga, \gw), - \infty
\leq \ga < \omega \leq \infty$.  Suppose there exists $K < \infty$ such that
\begin{align*}
\sup_{M_i \times (\ga, \gw)} \brs{\Rm(g_i)}_{g_i} \leq K.
\end{align*}
Then there exists a subsequence $\{(M_{i_j}, g_{i_j}(t), p_{i_j}) \}$ and a one
parameter family of complete pointed metric spaces $(X, d(t), x)$ such
that for each $t \in (\ga, \gw)$, $\{(M_{i_j}, d_{i_j}, p_{i_j}) \}$ converges
to
$(X, d(t), x)$ in the sense of $C^{\infty}$ local submersions.  The local
lifted metrics $h_y(t)$ are solutions to FOCF.  Furthermore, if there exists a
constant $\gd > 0$ so that
\begin{align*}
\inj_{g_i(0)}(p_i) \geq \gd
\end{align*}
then the limit space $(X, d(t), x)$ is a smooth $n$-dimensional Riemannian
manifold, and the limiting metric is the $C^{\infty}$ limit of the metrics
$g_i(t)$.
\end{cor}

Going yet further, one can actually obtain convergence in the category of
Riemannian groupoids.  Compactness of solutions to Ricci flow on Riemannian
groupoids satisfying a uniform curvature bound was established by Lott
\cite{Lott}.  Again, this is less a statement about Ricci flow as it is a
statement about one parameter families of metrics satisfying certain bounds on
curvature and its derivatives.  Therefore our smoothing estimates can be
employed to carry his argument over to fourth order curvature flows.  The
relevant background and definitions will be given in section \ref{corollaries}.

\begin{thm} \label{groupoidcompactness} Let $\{(M^n_i, g_i(t), p_i)\}$ be a
sequence of complete pointed
solutions of FOCF, where $t \in (\ga, \gw), - \infty
\leq \ga < \omega \leq \infty$.  Suppose there exists $K < \infty$ such that
\begin{align*}
\sup_{M_i \times (\ga, \gw)} \brs{\Rm(g_i)}_{g_i} \leq K.
\end{align*}
Then there exists a subsequence $\{(M_{i_j}, g_{i_j}(t), p_{i_j}) \}$ such that
$g_{i_j}(t)$ converges smoothly to a solution to FOCF $g_{\infty}(t)$ on an
$n$-dimensional \'etale groupoid $(G_{\infty}, p_{\infty})$ defined on $(\ga,
\gw)$.
\end{thm}

\begin{rmk} One could phrase everything we do here in the language of
groupoids, but in the interests of clarity we have mainly focused on the notion
of convergence of local submersions, which is probably more familiar to most.
\end{rmk}

\begin{rmk} For instances of FOCF which are gradient flows of positive quadratic
scale-invariant functionals,
blowup limits of noncollapsed finite time singularities are automatically
critical for the given functional.  An interesting question prompted by this
theorem is whether, for such gradient flows, the limits given by Theorem
\ref{groupoidcompactness} are still critical points.
\end{rmk}

As another corollary to Theorem \ref{globalsmoothing}, we are able to give an
improvement of Theorem \ref{oldLTE}.

\begin{cor} \label{existencecor} 
Let $(M^n, g)$ be a compact Riemannian manifold.  The
solution to (\ref{L2flow}) with initial condition $g$ exists on a maximal time
interval $[0, T)$.  Furthermore, if $T < \infty$ then
\begin{align*}
\limsup_{t \to T} \brs{\Rm}_{g(t)} = \infty.
\end{align*}
\end{cor}

\begin{rmk}
Theorem \ref{oldLTE} allowed for the possibility that the solution could
collapse with bounded curvature in finite time.  This would be a very
unfortunate outcome, leaving little hope to understand the singularity with
blowup arguments.  Corollary \ref{existencecor} ensures curvature blowup at a
finite singular time.
\end{rmk}

\begin{rmk} Comparing Corollary \ref{existencecor} with Theorem \ref{oldLTE}
you might think we have shown that the Sobolev constant is always bounded, which
is not the case.  What we show is actually $C^0$ control over the metric in the
presence of a curvature bound, a nontrivial statement for solutions to FOCF.
\end{rmk}

\begin{rmk}
An
important open problem is to show a noncollapsing estimate on the scale of
curvature akin to Perelman's estimate for Ricci flow.  Such a general estimate
for solutions to FOCF which are gradient flows of $L^2$ curvature energies would
in particular imply their long time existence on surfaces and three-manifolds
using
a blowup argument and Corollary \ref{existencecor}.
\end{rmk} 

\begin{rmk} We point out here that the methods of this paper should adapt
to higher-order (even higher than fourth) geometric flows on manifolds, and
other situations where the maximum principle is not available.  
A common feature of these equations is the presence of Sobolev space smoothing
estimates, which require Sobolev constant hypotheses to convert into pointwise
bounds to yield long time existence obstructions.  The procedure employed in
this paper allows you to remove some of these hypotheses.  Here is the overall
outline of the strategy:
\begin{itemize}
\item {Prove a precompactness result assuming $C^k$ bounds on curvature.}
\item {By applying the $L^2$ smoothing estimates on the covering metric, use the
compactness result and a blowup argument to show parabolic $C^k$ derivative
estimates in the presence of a curvature bound.}
\item{Turn these pointwise derivative estimates around the obtain a superior
compactness result only assuming a curvature bound.}
\item{Derive long time existence in the presence of a curvature bound using
the pointwise smoothing estimates.}
\end{itemize}
\end{rmk}

\begin{rmk} An interesting open problem is to localize the estimates of Theorem
\ref{globalsmoothing}.  Here again the fact that a curvature bound does not
imply $C^0$
control of the metric makes things particularly difficult.  This problem of
course does not arise for Ricci flow, and in fact, the blowup argument of
Theorem \ref{globalsmoothing} can be applied
locally to solutions of Ricci flow to get an estimate with particularly
clean dependencies on time, curvature, and distance to the boundary.
\end{rmk}

Our final theorem concerns nonsingular solutions to the $L^2$ gradient flow.  We
first need to introduce the volume normalized version of (\ref{L2flow}).  It
follows from (\ref{gradF}) that if $(M^n, g(t))$ is a solution to (\ref{L2flow})
then
\begin{align*}
\dt \Vol(g(t)) =&\ \frac{4 - n}{4} \mathcal F(g(t)).
\end{align*}
In particular, the initial value problem
\begin{gather} \label{vnL2flow}
\begin{split}
\dt g =&\ - \grad \FF + \frac{n-4}{2n} \frac{\mathcal F(g)}{\Vol(g)} \cdot
g\\
g(0) =&\ g_0,
\end{split}
\end{gather}
preserves the volume of the time dependent metrics.  One can check that for an
initial metric $g_0$, the corresponding solutions to (\ref{L2flow}) and
(\ref{vnL2flow}) differ by a rescaling in space and time.  Furthermore, equation (\ref{vnL2flow}) is, once suitably normalized, the gradient flow of the functional
\begin{gather} \label{Ftildef}
 \til{\FF}(g) = \Vol(g)^{\frac{4-n}{n}} \FF(g)
\end{gather}

\begin{defn} A solution $(M^n, g(t))$ to (\ref{vnL2flow}) is \emph{nonsingular}
if
it exists on $[0, \infty)$ with a
uniform bound on the curvature tensor.
\end{defn}

\begin{thm} \label{nonsingularthm} Let $(M^n, g(t))$ be a nonsingular solution
to (\ref{vnL2flow}).  Then either
\begin{itemize}
\item{ For all $p \in M$, $\limsup_{t \to \infty} \inj_p g(t) = 0$.}
\item{ There exists a sequence of times $t_i \to \infty$ such that $\{g(t_i) \}$
converges to a smooth metric on $M$ which is critical for $\til{\FF}$.}
\item{ There exists a sequence of points $(p_i, t_i), t_i \to \infty$ such that
$\{(M, g(t_i), p_i) \}$
converges to a complete noncompact finite volume metric which is critical for
$\til{\FF}$.}
\end{itemize}
\end{thm}

\begin{rmk} The same theorem holds for solutions to any FOCF which is the
gradient flow of a positive scale invariant functional, in particular the Calabi
flow (replacing the word
\emph{critical} with \emph{extremal} everywhere), and the flows of Bour
mentioned above.
\end{rmk}

Here is an outline of the rest of the paper.  In section \ref{wkconvsec} we
collect a number of background results on Gromov-Hausdorff convergence from
\cite{Glickenstein}.  Section \ref{proofwkwk} has the proof of Theorem
\ref{wkwkcompactness}.  In section \ref{localsobsmoothingsec} we prove local
smoothing estimates for solutions to FOCF in the presence of certain pointwise
bounds on curvature and the time derivative of the metric.  In section
\ref{localsmoothingsec} we give the proof of Theorem \ref{globalsmoothing}.  We
end in section \ref{corollaries} with the proofs of
Corollary \ref{wkcompactness}, \ref{existencecor}, and Theorem
\ref{nonsingularthm}.

\noindent \textbf{Acknowledgements:} The author would like to thank Aaron Naber
for some helpful discussions, and an anonymous referee for a very thorough
reading of a previous version of this paper.

\section{Weak Convergence} \label{wkconvsec}

In this section we recall certain definitions of weak convergence of Riemannian
manifolds.  Further background on these definitions may be found in
\cite{Fukaya}, \cite{Glickenstein}.  

\begin{defn} A topological space $G$ is a \emph{pseudogroup} if there exist
pairs $(a, b) \in G \times G$ such that a product $ab \in G$ is defined
satisfying
\begin{enumerate}
\item {If $ab, bc, (ab)c, a(bc)$ are all defined, then $(ab)c = a(bc)$.}
\item{ If $ab$ is defined then for every neighborhood $W$ of $ab$ there are
neighborhoods $a \in U$ and $b \in V$ such that for all $x \in U, y \in V$, $xy$
is defined and $xy \in W$.}
\item{ There is an element $e \in G$ such that for all $a \in G$, $ae$ is
defined and $ae = a$.}
\item{ If for $(a, b) \in G \times G$, $ab$ is defined and $ab = e$, then $a$ is
a left-inverse for $b$ and we write $a = b^{-1}$.  If $b$ has a left inverse,
then for every neighborhood $U$ of $b^{-1}$ there is a neighborhood $V$ of $b$
such that every $y \in V$ has a left inverse $y^{-1} \in U$.}
\end{enumerate}
\end{defn}

\noindent The canonical example of a pseudogroup is a Lie group germ.

\begin{defn} A pseudogroup $G$ is a \emph{Lie group germ} if a neighborhood of
the identity of $G$ can be given a differentiable structure such that the
operations of group multiplication and inversion are differentiable maps when
defined.
\end{defn}

\begin{defn} \label{wkconvdef} Fix $k \in \mathbb R \cup \{\infty \}$.  A
sequence of pointed
$n$-dimensional Riemannian manifolds $\{(M_i^n, g_i, p_i)\}$ \emph{locally
converges to a pointed metric space $(X, d, x)$ in the sense of $C^{k}$-local
submersions at x} if there is a Riemannian metric $h$ on an open neighborhood $0
\in V \subset \mathbb R^n$, a pseudogroup $\Gamma$ of local isometries of $(V,
h)$ such that the quotient is well-defined, an open set $U \subset X$ and maps
\begin{align*}
\phi_i : (V, 0) \to (M_i, p_i)
\end{align*}
satisfying
\begin{enumerate}
\item { $\{(M_i, d_{g_i}, p_i) \}$ converges to $(X, d, x)$ in the pointed
Gromov-Hausdorff topology,}
\item{the identity component of $\Gamma$ is a Lie group germ,}
\item{$V / \Gamma, \bar{d}_h \cong (U, d)$ where $\bar{d}_h$ is the induced
distance function on the quotient,}
\item{$(\phi_i)_*$ is nonsingular on $V$ for all $i$,}
\item{$h$ is the $C^{k}$ limit of $\phi_i^* g_i$, in the sense of uniform
convergence on compact sets of the first $k$ derivatives.  As $k \in \mathbb R$
this is meant in the usual H\"older sense.}
\end{enumerate}
\end{defn}

\begin{defn} Fix $k \in \mathbb N \cup \{\infty \}$.  A sequence of pointed
$n$-dimensional Riemannian manifolds $\{(M_i, g_i, p_i \}$ \emph{converges to a
pointed metric space $(X, d, x)$ in the sense of $C^{k}$ local submersions} if
for every $y \in X$ there are points $q_i \in M_i$ such that $\{(M_i, g_i,
q_i)\}$ locally converges to $(X, d, y)$ in the sense of $C^{k}$ local
submersions.
\end{defn}

We now recall the compactness theorem of Glickenstein, which we will adapt in
the next section to one-parameter families of metrics satisfying certain bounds
on derivatives of curvature.

\begin{thm} \label{glick1} (\cite{Glickenstein} Theorem 3) Let $C_k > 0$ be
constants for $k \in \mathbb N$ and $\{(M_i^n, g_i(t), p_i) \}_{i =
1}^{\infty}$, where $t \in [0, T]$, be a sequence of pointed solutions to the
Ricci flow on complete manifolds such that
\begin{align*}
\brs{\Rm}_{g_i(t)} \leq&\ 1\\
\brs{\N_i^k \Rm}_{g_i(t)} \leq&\ C_k,
\end{align*}
for all $i, k \in \mathbb N$ and $t \in [0, T]$.  Then there is a subsequence
$\{(M_{i_k}, g_{i_k}(t), p_{i_k})\}_{k = 1}^{\infty}$ and a one parameter family
of complete pointed metric spaces $(X, d(t), x)$ such that for each $t \in [0,
T]$, $(M_{i_k}, d_{g_{i_k}}(t), p_{i_k})$ converges to $(X, d(t), x)$ in the
sense of $C^{\infty}$-local submersions and the metrics $h_y(t)$ are solutions
to the Ricci flow equation.
\end{thm}

A fundamental result at the root of the proof of Theorem \ref{glick1} and our
adaptation Theorem \ref{wkwkcompactness} below is the next theorem on
convergence of families of Riemannian manifolds.  

\begin{thm} \label{Glickcomp} (\cite{Glickenstein} Theorem 20, Proposition 21) 
Let $\{(M_i, g_i(t), p_i)) \}$, $t \in [0, T]$ be a sequence of pointed
Riemannian manifolds of dimension $n$ satisfying that for every $\gd > 0$ there
is $\ge > 0$ such that if $s, t \in [0, T]$, $\brs{s - t} < \ge$ then
\begin{align*}
(1 + \gd)^{-1} g_i(s) \leq g_i(t) \leq (1 + \gd) g_i(t_0)
\end{align*}
for all $i > 0$, and satisfying
\begin{align*}
\Rc(g_i(t)) \geq - A g_i(t).
\end{align*}
Then there is a subsequence $\{(M_{i_k}, g_{i_k}(t), p_{i_k}) \}$ and a
$1$-parameter family of complete pointed metric spaces $(X_{\infty}(t),
d_{\infty}(t), x_{\infty}(t))$ such that for every $t \in [0, T]$ the
subsequence converges to $(X_{\infty}(t), d_{\infty}(t), x_{\infty}(t))$ in the
Gromov-Hausdorff topology.  If furthermore
\begin{align*}
\brs{\dt g} \leq C
\end{align*}
for all $t \in [0, T]$ then $(X_{\infty}(t), x_{\infty}(t))$ is homeomorphic to
$(X_{\infty}(0), x_{\infty}(0))$ for all $t \in [0, T]$.
\end{thm}

\section{A weak compactness result} \label{proofwkwk}

The purpose of this section is to prove Theorem \ref{wkwkcompactness}, a weak
compactness theorem for one-parameter families of Riemannian metrics satisfying
certain curvature bounds.  As mentioned in the introduction, the proof is a
direct adaptation of the main result of \cite{Glickenstein}, therefore we only
indicate the steps required to modify that proof and refer the reader to
\cite{Glickenstein} for a more thorough discussion.

\begin{thm} \label{wkwkcompactness} Fix $k \in \mathbb \mathbb N \cup \{\infty
\}, k \geq 3$.  Let $\{(M_i, g_i(t), p_i)\}$ be a sequence of complete pointed
solutions of FOCF on manifolds $M_i$, where $t \in (\ga, \gw), - \infty \leq \ga
< \omega \leq \infty$.  Suppose that for all $0 \leq l \leq k$ there exists a
constant $C_l > 0$ such that
\begin{align*}
\brs{\N_i^l \Rm(g_i)}_{g_i} \leq C_l.
\end{align*}
There exists a subsequence $\{(M_{i_j}, g_{i_j}(t), p_{i_j}) \}$ which converges
to a one parameter family of complete pointed metric spaces $(X, d(t), x)$ such
that for each $t \in (\ga, \gw)$, $\{(M_{i_j}, d_{i_j}, p_{i_j}) \}$ converges
to
$(X, d(t), x)$ in the sense of $C^{k-2-\gb}$ local submersions.  The local
lifted
metrics $h_y(t)$ are solutions to FOCF.
\end{thm}

\begin{rmk} A key point of the statement is that the
local submersion structure, i.e. the local groups of isometries and quotient
maps of Definition \ref{wkconvdef}, are independent of time.
\end{rmk}

\begin{lemma} \label{gcontinuity} Suppose $(M^n, g(t))$ is a one parameter
family of metrics and suppose
\begin{align*}
\sup_{M \times [0, T)} \brs{\dt g}_{g(t)} \leq A.
\end{align*}
Then for all $s, t \in [0, T)$,
\begin{align*}
e^{- A \brs{t - s}} g(s) \leq g(t) \leq e^{A \brs{t - s}} g(s).
\end{align*}
\begin{proof}  This is a straightforward estimate, see for instance \cite{Chow}
Lemma 6.49.
\end{proof}
\end{lemma}

\begin{lemma} \label{ballgrowth}  Suppose $(M^n, g(t))$ is a one parameter
family of metrics and suppose
\begin{align*}
\sup_{M \times [0, T)} \brs{\dt g}_{g(t)} \leq A.
\end{align*}
Then for all $\rho > 0$, 
\begin{align*}
B_{g(t)}(0, r_A(t) \rho) \subset&\ B_{g(0)}(0, \rho),\\
B_{(g(0)}(0, r_A(t) \rho) \subset&\ B_{g(t)}(0, \rho),
\end{align*}
where
\begin{align*}
r_A(t) = \frac{1}{1 + (e^{At} - 1)^{\frac{1}{2}}}.
\end{align*}
\begin{proof} Again this is a straightforward estimate, and one can consult
\cite{Glickenstein} Proposition 19 for the argument for Ricci flow.
\end{proof}
\end{lemma}

\begin{lemma} \label{curvaturetometricbnds} Let $M$ be a Riemannian manifold
with metric $g$.  Let $K$ denote a compact subset of $M$, and $\{g_i \}$ a
sequence of solutions to FOCF defined on open neighborhoods of $K \times [\ga,
\gb ]$, where $0 \in [\ga, \gb]$.  Let $D$, $D_i$ and $\brs{\cdot}$,
$\brs{\cdot}_i$ denote the Levi Civita connections and norms of $g$ and $g_i$
respectively.  Fix $N \geq 0$, and suppose
\begin{enumerate}
\item{There is a constant $C > 0$ so that on $K$ one has
\begin{align*}
\frac{1}{C} g(0) \leq g_i(0) \leq C g(0).
\end{align*}}
\item{For $j \leq N$, the covariant derivatives of $\{g_i \}$ with respect to
$D$
are uniformly bounded at $t = 0$ on $K$, i.e.
\begin{align*}
\brs{D^j g_i} \leq C_j.
\end{align*}}
\item{For $j \leq N + 2$, the $j$-th covariant derivative of $\{\Rm_i\}$ is
bounded
with respect to $\{g_i\}$ on $K \times [\ga, \gb]$, i.e.
\begin{align*}
\brs{D_i^j \Rm_i}_i \leq C_j.
\end{align*}}
\end{enumerate}
Then the metrics $\{g_i\}$ are uniformly bounded with respect to $g$ on $K
\times [\ga, \gb]$, and moreover for $j \leq N$ the $j$-th covariant derivative
of $\{g_i\}$ with respect to $D$ is uniformly bounded on $K \times [\ga, \gb]$,
i.e.
\begin{align*}
\brs{D^j g_i} \leq \til{C}_j,
\end{align*}
where these bounds all depend only on the assumed bounds and the dimension.
\begin{proof} This is \cite{Hamilton} Lemma 2.4, adapted to the case of FOCF,
and also to the case where one is only assuming a finite number of derivatives
of curvature are bounded.  The proof relies on a straightforward integration in
time using the derivative bounds.
\end{proof}
\end{lemma}

\noindent We now proceed with the proof of Theorem \ref{wkwkcompactness}, which
closely follows the proof of \cite{Glickenstein} Theorem 3.

\begin{proof}
Recall that the statement assumed uniform pointwise bounds on the curvature and
at least its first $3$ derivatives.  In particular the Ricci curvature is
bounded
below and $\brs{\dt g}$ is bounded.  By applying Lemma \ref{gcontinuity} we see
that the hypotheses of Theorem \ref{Glickcomp} are satisfied, and thus there is
a one-parameter family of complete pointed metric spaces 
$(X, d(t), x)$ such that $\{(M_i, g_i(t), p_i)\}$ converges in the pointed
Gromov-Hausdorff topology to $(X, d(t), x)$.  It remains to understand
the local
submersion structure of this limit space $(X, d(t), x)$, and the proof
is at this point a direct adaptation of the main result of \cite{Glickenstein},
the only difference being convergence in $C^{k-2-\ga}$, and so we refer the
reader there for details.
\end{proof}

\section{Local Smoothing Estimates in Sobolev Spaces}
\label{localsobsmoothingsec}

In this section we derive local estimates for the $H_2^l$ norms of curvature for
a solution to (\ref{generalflow}).  We will compute some evolution equations,
but first let we introduce a helpful piece of notation.  Given a tensor $A$, we
we will write  $P_s^m(A)$ for any universal expression of the form
\begin{align*}
P_s^m(A) = \sum_{i_1 + \dots + i_s = m} \N^{i_1} A * \dots * \N^{i_s} A.
\end{align*}

\begin{lemma} Let $(M^n, g(t))$ be a solution to FOCF.  Then
\begin{align*}
\dt \N^k \Rm =&\ - \gD^2 \N^k \Rm + P_2^{k+2}(\Rm) + P_3^k (\Rm).
\end{align*}
\begin{proof} The proof is identical to \cite{StreetsL21} Proposition 4.3.
\end{proof}
\end{lemma}

For now we fix a given solution $(M^n, g)$ to FOCF, and fix some further data. 
In particular, fix $\gg \in C_0^{\infty}(M)$.  Suppose that $K$ is a constant
such
that
\begin{align} \label{localbound}
\sup_{\supp \gg \times [0, T]} \left\{ \brs{\Rm} + \brs{\dt g} + \brs{\N \dt g}
\right\} \leq&\ K.
\end{align}

\begin{lemma} \label{cutoffboundlemma} Let $(M^n, g(t))$ denote a one-parameter
family of Riemannian
metrics on $[0, T]$.  Let $\gamma \in C_0^{\infty}(M)$.  There are constants
$C_1$ and $C_2$ such that
\begin{align*}
\brs{d \gamma}_{g(t)} \leq&\ C_1 \left(\brs{\dt g}, T, \brs{d \gamma}_{g(0)}
\right),\\
\brs{\N \N \gamma}_{g(t)} \leq&\ C_2 \left(\brs{\dt g},
\brs{\N \dt g}, T, \brs{d \gamma}_{g(0)}, \brs{\N^2 \gamma}_{g(0)} \right).
\end{align*}
\begin{proof} This first estimate is straightforward.  For the second we express
the
time-dependent Hessian as
\begin{align*}
\left(\N \N \gg\right)_{ij} =&\ \del_i \del_j \gg - \gG_{ij}^k \del_k \gg.
\end{align*}
Differentiating with respect to time yields
\begin{align*}
\dt \brs{\N \N \gg}^2_{g(t)} =&\ \N \dt g * d \gg * \N \N \gg + \dt g * \N^2
\gg^{*2}\\
\leq&\ C \left( \brs{d \gg}, \brs{\N^2 \gg}, \brs{\dt g}, \brs{\N \dt g}
\right).
\end{align*}
The result follows.
\end{proof}
\end{lemma}

\noindent In our setup the above lemma implies the estimates
\begin{gather} \label{cutoffestimates}
\begin{split}
\brs{d \gg}_{g(t)} \leq&\ C \left(K, T, \brs{d \gg}_{g(0)} \right),\\
\brs{\N \N \gg}_{g(t)} \leq&\ C \left(K, T, \brs{d \gg}_{g(0)}, \brs{\N^2
\gg}_{g(0)} \right).
\end{split}
\end{gather}

\noindent We now recall an estimate from \cite{StreetsL21}.  There is a typo we
correct here in
where line (4.7) in that paper only assumes a bound on the cutoff function for
the
initial metric.  Also, we are applying the result to any
solution to FOCF, not just solutions to (\ref{L2flow}), but the proof is
identical.

\begin{prop} (\cite{StreetsL21} Corollary 5.2) \label{L2estimate} Let $(M^n,
g(t))$ be a solution
to FOCF on $[0, T]$.  Fix $\gg \in C_0^{\infty}(M)$ and suppose $K$ and $L$ are
constants such that
\begin{align*}
\sup_{\supp \gg \times [0, T]} \brs{\Rm} \leq K,\\
\sup_{[0, T]} \left(\brs{d \gg} + \brs{\N \N \gg} \right) \leq L.
\end{align*}
Fix $W = \N^m \Rm$.  Then for any $s \geq 2m + 4$, there exists a constant $C =
C(m, s, L)$ such that
\begin{align*}
\dt \int_M \brs{W}^2 \gg^s + \frac{1}{2} \int_M \brs{\N^2 W}^2 \gg^s \leq C K^2
\int_M \brs{W}^2 \gg^s + C \left(1 + K^2 \right) \brs{\brs{\Rm}}_{L^2(\{\supp
\gg \})}^2. 
\end{align*}
\end{prop}

Using this estimate we can prove a local smoothing estimate in Sobolev space
norms which we will employ in the proof of Theorem \ref{globalsmoothing}.

\begin{thm} \label{localSobestimates} Let $(M^n, g(t))$ be a solution to FOCF on
$[0, T]$.  Fix $r > 0$ and suppose $x \in M$ satisfies
\begin{align} \label{spacetimebound}
\sup_{[0, T] \times B_{g(T)}(x, r)} \left \{ \brs{\Rm} + \brs{\dt g} + \brs{\N
\dt g} \right\} \leq
K.
\end{align}
Given $m \geq 0$ there exists $C = C\left(m, \frac{1}{r}, T, K
\right)$ such
that
\begin{align*}
\brs{\brs{\N^m \Rm_{g(t)}}}_{L^2\left(B_{g(T)}\left(x, r \right) \right)}^2
\leq \frac{C}{t^m} \sup_{[0, T]} \brs{\brs{\Rm}}_{L^2({B_{g(T)}(x, 2 r)})}^2.
\end{align*}
\begin{proof} We show the result for $m = 2n \geq 0$ even, and the result
follows for $m$ odd by interpolation.  Fix $m \geq 0$ even, and let $\gb_k$ be
constants to be determined.  Let $\gg$ denote a cutoff function for the ball of
radius $r$ with respect to the metric $g(T)$.  This $\gg$ satisfies
\begin{align*}
\brs{d \gg}_{g(T)} + \brs{\N \N \gg}_{g(T)} \leq C \left(\frac{1}{r} \right).
\end{align*}
Furthermore, by (\ref{spacetimebound}) and Lemma \ref{cutoffboundlemma} we
conclude that
\begin{align*}
\sup_{[0, T]} \left( \brs{d \gg} + \brs{\N \N \gg} \right) \leq L \left(T, K,
\frac{1}{r}
\right).
\end{align*}
Now let
\begin{align*}
\Phi := \sum_{k = 0}^n \gb_k t^{2k} \brs{\brs{ \gg^{m + 2} \N^{2k} \Rm
}}_{L^2}^2.
\end{align*}
for constants $\gb_k$ to be determined, with $\gb_n = 1$.
It follows from Proposition \ref{L2estimate} that 
\begin{align*}
\frac{d}{dt} \Phi \leq&\ \sum_{k = 1}^{n} \brs{\brs{\gg^{m+2} \N^{2k}
\Rm}}_{L^2}^2 \left(C K^2 \gb_{k} t^{2k} - \frac{1}{2} \gb_{k-1} t^{2k-2}
+ 2 k \gb_k t^{2k-1} \right)\\
&\ + C(m, s, L, K, T, \gb_i) \brs{\brs{\Rm}}_{L^2(\{\supp \gg\})}^2.
\end{align*}
It is clear that by an appropriate inductive choice of the constants $\gb_i$
with respect to the constants $K, L, T$ and $m$ we obtain
\begin{align*}
\dt \Phi \leq&\ C(m, s, L, K, T) \brs{\brs{\Rm}}_{L^2(\{\supp \gg \})}^2.
\end{align*}
Integrating this ODE yields the result for even $m$, and for odd $m$ the result
follows using an interpolation inequality.
\end{proof}
\end{thm}

\section{Proof of Theorem \ref{globalsmoothing}} \label{localsmoothingsec}

Here is a sketch of the proof of Theorem \ref{globalsmoothing}.  Suppose there
is no such constant
$C$.  Then one has a sequence of solutions $\{(M_i, g_i(t)\}$
(note that the topology of $M_i$ is unknown) to FOCF with bounded
curvature such that there exist points $(x_i, t_i)$ violating 
(\ref{globalsmoothingestimate}).  Precisely by the scale-invariant nature of the
claimed
estimates, one can blow up at the scale of $\N^m \Rm$ and argue to get a
well-defined
local limit converging in the sense of $C^{\infty}$ submersions by the weak
compactness theorem.  This limit has vanishing curvature by construction.  But
also by construction the $m$-th covariant derivative is nonzero at the center
point, a contradiction.

\begin{proof}  Fix $m, n > 0$, and define the function
\begin{align*}
f_m(x, t, g) := \sum_{j = 1}^m \brs{\N^j \Rm}_{g(t)}^{\frac{2}{2 + j}}(x).
\end{align*}
We will show that given a complete solution to FOCF as in the hypotheses,
for all sufficiently large $m$ there is a constant $C = C(m, n)$ such that
\begin{align} \label{gtechnicalsmoothing}
f_m(x, t, g) \leq C \left(K + \frac{1}{t^{\frac{1}{2}}} \right).
\end{align}
This clearly suffices to prove the theorem.  It suffices for the proof here to
assume $m \geq 3$.  Suppose the claim was false.  Take then a sequence $ \{(M_i,
g_i(t)) \}$ of complete solutions to FOCF satisfying the hypotheses of the
theorem,
together with points $(x_i, t_i)$ satisfying
\begin{align*}
\lim_{i \to \infty} \frac{f_m(x_i, t_i, g_i)}{K + \frac{1}{t_i^{\frac{1}{2}}}} =
\infty,
\end{align*}
Since each solution is smooth it follows that $t_i > 0$.  Without loss of
generality we may choose the
points $(x_i, t_i)$ such
that
\begin{align} \label{gsupcondition}
\frac{f_m(x_i, t_i, g_i)}{K + \frac{1}{t_i^{\frac{1}{2}}}} =
\sup_{M_i \times \left(0, T \right] } \frac{f_m(x, t, g_i)}{K +
\frac{1}{t^{\frac{1}{2}}}}.
\end{align}
Let $\gl_i = f_m(x_i, t_i)$, and set
\begin{align*}
\til{g}_i = \gl_i g\left(t_i + \frac{t}{\gl_i^2} \right).
\end{align*}
Let us make some observations about $\til{g}_i$ which make clear why the
estimates of the theorem take the form they do and moreover why we have made the
choices above.  First, observe that by construction the solution $\til{g}_i$
exists on the time interval $[- t_i \gl_i^2, 0]$.  But since
\begin{align*}
t_i^{\frac{1}{2}} \gl_i = \frac{f_m(x_i, t_i)}{\frac{1}{t_i^{\frac{1}{2}}}} \geq
\frac{f_m(x_i, t_i)}{K + \frac{1}{t_i^{\frac{1}{2}}}},
\end{align*}
and the right hand side above goes to infinity as $i \to \infty$
we conclude that the solutions $\til{g}_i$ exist on $[-1, 0]$ for all
sufficiently large $i$.  Next observe that 
\begin{align*}
\frac{\gl_i}{K} = \frac{f_m(x_i, t_i)}{K} \geq \frac{f_m(x_i, t_i)}{K +
\frac{1}{t_i^{\frac{1}{2}}}}
\end{align*}
and again the right hand side goes to infinity as $i \to \infty$.  It follows
that 
\begin{align} \label{gsmoothingloc10}
\lim_{i \to \infty} \brs{\til{\Rm_i}} \leq \lim_{i \to \infty} \frac{K}{\gl_i} =
0.
\end{align}
\noindent By construction clearly
\begin{align} \label{gfnonzero}
\til{f}_m(x_i, 0) = 1.
\end{align}
Also, we observe that for $(x'_i, t'_i) \in M_i \times [-1, 0]$, using
(\ref{gsupcondition}) one has
\begin{gather} \label{gfbound}
\begin{split}
\til{f}_m(x_i', t_i') =&\ \frac{ f_m \left(x_i', t_i + \frac{t_i'}{\gl_i^2}
\right)}{\gl_i}\\
=&\ \frac{ f_m \left(x_i', t_i + \frac{t_i'}{\gl_i^2} \right)}{f_m(x_i, t_i)}\\
\leq&\ \frac{K + \left(t_i + \frac{t_i'}{\gl_i^2} \right)^{-\frac{1}{2}}}{K +
t_i^{-\frac{1}{2}}}\\
\leq&\ 1.
\end{split}
\end{gather}

Thus the sequence of pointed solutions $\{(M_i, \til{g}_i, x_i) \}$ has a
uniform
estimate on the first $m \geq 3$ derivatives of curvature on $M_i \times [-1,
0]$.  By Theorem \ref{wkwkcompactness} we have a one-parameter family of
pointed metric spaces $(X, d(t), x)$ and a subsequence (denoted with the same
index) such that $\{(M_i, \til{g}_i(t), x_i) \}$
converges to $(X, d(t), x)$ in
the sense of $C^{m-2-\ga}$-local submersions. Actually, for our purposes here,
the global convergence statement of Theorem \ref{wkwkcompactness} is not needed,
only the pointwise statement that the pullbacks of $g_i(t)$ by the exponential
map at $x_i$ converge on a ball in $\mathbb R^n$.  In particular there is a
sequence
$\{\til{h_i}\}$
of one-parameter families of local liftings of $\til{g}_i(t)$ near $x_i$,
defined on some ball $B(0, r) \subset \mathbb R^n$, converging to a
one-parameter family $h_{\infty}(t)$, as guaranteed by
Theorem \ref{wkwkcompactness}.  So far the convergence to
$\til{h}_{\infty}$ is only in the $C^{m - 2 - \ga}$ topology, but we can improve
this to $C^{\infty}$ convergence using Theorem \ref{localSobestimates}.  As the
metrics $\til{h}_i$ are defined using the exponential map at $x_i$, and the
curvature is uniformly bounded, we have that the metrics $\til{h}_i(0)$ are
uniformly $C^0$ equivalent to the Euclidean metric on $B(0, r)$.
It follows that the Sobolev constant of
$B_{\til{h}_i(0)}(0, r)$ is uniformly bounded.  Also, since $m \geq 3$, by
(\ref{gfbound}) we have a uniform
bound on $\dt g$ and $\N \dt g$ on $[-1, 0]$ and so it follows from Theorem
\ref{localSobestimates} that the $H_2^p$ norms of the curvature of
$\til{h}_i(s), s \geq - \frac{1}{2}$ are uniformly bounded for all $p > 0$. 
Using the Sobolev constant bound of $\til{h}_i(0)$, it follows that the $C^k$
norms of the
curvature of $\til{h}_i(0)$ are uniformly bounded for all $k$, and thus by
taking a further
subsequence we conclude $C^{\infty}$ convergence to $\til{h}_{\infty}(0)$.  We
conclude from this convergence and (\ref{gfnonzero}) that
\begin{align*}
f_m(0, 0, \til{h}_{\infty}) = 1.
\end{align*}
However, we conclude from (\ref{gsmoothingloc10}) that on $B(0, r)$,
\begin{align*}
\Rm_{\til{h}_{\infty}} \equiv 0.
\end{align*}
This is a clear contradiction and so the theorem follows.
\end{proof}

\section{Compactness Theorems}
\label{corollaries}

\noindent We begin with the proof of Corollary \ref{wkcompactness}.

\begin{proof} By a diagonalization argument it suffices to find a convergent
subsequence on a sequence of closed intervals  whose endpoints approach $\ga$
and $\gw$.  By Theorem \ref{globalsmoothing}, the uniform curvature estimate
implies a uniform bound on all derivatives of curvature for any closed interval
contained inside $(\ga, \gw)$.  Theorem \ref{wkwkcompactness} then guarantees
the existence of a subsequence converging in the sense of $C^{\infty}$-local
submersions.  Finally, if the injectivity radii at $p_i$ are uniformly bounded
below, we claim the limit is a smooth manifold of the same dimension.  At this
point we have all of the estimates required in Hamilton's proof of the
compactness of Ricci flow solutions under the hypotheses of bounded curvature
and injectivity radius, and so the proof there carries over with trivial
changes.  Specifically, one can apply \cite{Hamilton} Theorem 2.3 to our given
sequence to conclude the existence of a subsequence where the Riemannian
manifolds $\{M_i, g_i(0), p_i\}$ converge to $\{M_{\infty}, g_{\infty}, 0\}$. 
Next one can apply Lemma \ref{curvaturetometricbnds} (an easy adaptation of
\cite{Hamilton} Lemma 2.4) to conclude that there is a further subsequence which
converges uniformly on compact subsets of $M \times (\ga, \gw)$.
\end{proof}

Now we give the proof of Theorem \ref{groupoidcompactness}.  We start with the
main definition of a smooth \'etale groupoid.  We do not give the more general
definitions.  One may consult \cite{Lott} and the references therein for further
background information.

\begin{defn} A \emph{smooth \'etale groupoid} is a pair of smooth manifolds
$G^{(0)}$, $G^{(1)}$ together with 
\begin{enumerate}
\item {a smooth embedding
\begin{align*}
e : G^{(0)} \to G^{(1)},
\end{align*}}
\item{ Source and range local diffeomorphisms $s, r : G^{(1)} \to G^{(0)}$
satisfying
\begin{align*}
s \circ e = r \circ e = \Id_{G^{(0)}}.
\end{align*}
}
\item{ A multiplication $G^{(1)} \times G^{(1)} \to G^{(1)}$ such that
\begin{enumerate}
\item{ $\gg_1 \gg_2$ is defined if and only if $s(\gg_1) = r(\gg_2)$, in which
case $s(\gg_1 \gg_2) = s(\gg_2)$ and $r(\gg_1 \gg_2) = r(\gg_1)$.}
\item{ $(\gg_1 \gg_2) \gg_3 = \gg_1(\gg_2 \gg_3)$ whenever both sides are well
defined.}
\item{ $\gg s(\gg) = r(\gg) \gg = \gg$.}
\end{enumerate}
}
\end{enumerate}
\end{defn}
Since we are only concerned here with smooth \'etale groupoids we will simply
refer to these as \emph{groupoids} for the remainder of this paper.
Given a groupoid, the \emph{orbit of $x$} is
\begin{align*}
\mathcal O_x = s(r^{-1}(x)).
\end{align*}
A \emph{pointed groupoid} is a groupoid together with a marked orbit $\mathcal
O_x$.  Associated to a groupoid is a pseudogroup of diffeomorphisms of
$G^{(0)}$. Specifically, given $\gg \in G^{(1)}$, there exists a neighborhood
$U$ of $\gg$ in $G^{(1)}$ such that $\{(s(\gg'), r(\gg') | \gg' \in U \}$ is the
graph of a diffeomorphism between neighborhoods of the source and range of
$\gg$.  The pseudogroup $\mathcal P$ associated to a groupoid is generated by
such local diffeomorphisms.  A groupoid is \emph{Riemannian} if there is a
Riemannian metric $g$ on $G^{(0)}$ such that the elements of $\mathcal P$ act as
Riemannian isometries.  Finally, the \emph{dimension} of a groupoid will be the
dimension of the smooth manifold $G^{(0)}$.

With these preliminaries in place we can define the notion of convergence of
Riemannian groupoids.

\begin{defn} Let $\{(G_i, \mathcal O_{x_i}) \}$ be a sequence of pointed
$n$-dimensional Riemannian groupoids, and let $(G_{\infty}, \mathcal
O_{x_{\infty}})$ be another pointed Riemannian groupoid.  Let $J_{\infty}$ be
the groupoid of $1$-jets of local diffeomorphisms of $G_{\infty}^{(0)}$.  We say
that \emph{$\{(G_i, \mathcal O_{x_i}) \}$ converges to $(G_{\infty}, \mathcal
O_{x_{\infty}})$ in the pointed smooth topology} if for all $R > 0$,
\begin{enumerate}
\item {There are pointed diffeomorphisms $\phi_{i, R} : B_R(\mathcal
O_{x_{\infty}}) \to B_R(\mathcal O_{x_i})$ defined for sufficiently large $i$ so
that
\begin{align*}
\lim_{i \to \infty} \phi_{i,R}^* g_i |_{B_R(\mathcal O_{x_i})} =
g_{\infty}|_{B_R(\mathcal O_{x_{\infty}})}.
\end{align*}}
\item{After applying $\phi_{i, R}$, the images of
$s_i^{-1}(\bar{B_{\frac{R}{2}}(\mathcal O_{x_i})}) \cap
r_i^{-1}(\bar{B_{\frac{R}{2}}(\mathcal O_{x_i})})$ in $J_{\infty}$ converge in
the Gromov-Hausdorff topology to the image of
$s_{\infty}^{-1}(\bar{B_{\frac{R}{2}}(\mathcal O_{x_{\infty}})}) \cap
r_{\infty}^{-1}(\bar{B_{\frac{R}{2}}(\mathcal O_{x_{\infty}})})$}
\end{enumerate}
\end{defn}

\noindent Note that the groupoid of $1$-jets of local diffeomorphisms of an
\'etale
groupoid is not itself \'etale, though this point will not concern us here. 
Also, note of course that we are identifying elements of $G^{(1)}_i$ first with
local diffeomorphisms of $G^{(0)}_i$, and then considering the elements of $J_i$
they generate.

We can now give the proof of Theorem \ref{groupoidcompactness}.
\begin{proof} The main content is \cite{Lott} Proposition 5.9:
\begin{prop} Let $\{(M_i, g_i, p_i) \}$ be a sequence of complete
$n$-dimensional Riemannian manifolds.  Suppose that for all $k, R \geq 0$ there
is a constant $C_{k, R}$ such that for all $i$,
\begin{align*}
\sup_{B_{R}(p_i)} \brs{\N^k \Rm(g_i)} \leq C_{k, R}.
\end{align*}
Then there is a subsequence of $\{(M_i, g_i, p_i) \}$ which converges in the
pointed smooth topology to an $n$-dimensional Riemannian groupoid $(G_{\infty},
\mathcal O_{\infty})$.
\end{prop}
With this we can apply the same strategy of Corollary \ref{wkcompactness} (i.e.
the strategy of \cite{Hamilton} Theorem 1.2).  Given the pointwise smoothing
estimates, the sequence of manifolds $\{(M_i, g_i(0), p_i) \}$ satisfies the
hypotheses of the above proposition.  Therefore a subsequence converges to some
pointed Riemannian groupoid $(G_{\infty}, g_{\infty}(0), \mathcal
O_{x_{\infty}})$.  One can repeat the argument of Lemma
\ref{curvaturetometricbnds} (\cite{Hamilton} Lemma 2.4) in the groupoid setting
(\cite{Lott} Theorem 5.12) to obtain a further subsequence which converges for
all times the flow exists. The FOCF equation of course passes the limit
solution.  
\end{proof}

\section{Curvature blowup and nonsingular solutions}

\noindent We begin with the proof of Corollary \ref{existencecor}
\begin{proof}  Suppose the maximal existence time of the solution
is $T < \infty$ but
\begin{align*}
\sup_{[0, T)} \brs{\Rm} < \infty.
\end{align*}
By Theorem \ref{globalsmoothing} one has
\begin{align*}
\sup_{[0, T)} \brs{\grad \FF} = C < \infty.
\end{align*}
It follows from Lemma \ref{gcontinuity} that 
\begin{align*}
e^{-C t} g(0) \leq g(t) \leq e^{C t} g(0).
\end{align*}
In particular, this $C^0$ equivalence of metrics clearly implies that
\begin{align*}
\lim_{t \to T} C_S(g(t)) < \infty
\end{align*}
By Theorem \ref{oldLTE} we conclude that $T$ cannot be the maximal existence
time.  This contradicts the hypothesis that the curvature was bounded, therefore
the corollary follows.
\end{proof}

\noindent We end with the proof of Theorem \ref{nonsingularthm}.

\begin{proof} First observe that since the functional $\til{\FF}$ is bounded
below by zero and nonincreasing along a solution to (\ref{vnL2flow}) we have that
\begin{align} \label{flowtocritical}
\int_{0}^{\infty} \int_M \brs{\grad \til{\FF}}^2 dV dt =&\ \til{\FF}(g(0)) - \lim_{t
\to \infty} \til{\FF}(g(t)) < \infty.
\end{align}
It follows that we may choose a sequence of times $t_i \to \infty$ such that
\begin{align} \label{convtocrit}
\lim_{i \to \infty} \brs{\brs{\grad \til{\FF}(g(t_i))}}_{L^2} = 0.
\end{align}
Now let us assume without loss of generality that the first case does not occur.
 Then there exists $p \in M$ and $\gd > 0$ such that $\inj_{g(t_i)}(p) \geq \gd > 0$. 
It follows from Corollary \ref{wkcompactness} that there is a subsequence, also
denoted $t_i$, such that $\{(M, g(t_i + t), p) \}$ converges to a new solution
$(M_{\infty}, g_{\infty}(t), p_{\infty})$ to (\ref{vnL2flow}).  Since the
convergence is $C^{\infty}$ on compact sets, it follows from
(\ref{convtocrit}) that the limiting metric $g_{\infty}(t) = g_{\infty}$ is
critical for $\til{\FF}$.  If $M_{\infty}$ is compact, it follows that $M_{\infty}$ is
diffeomorphic to $M$ and thus the second alternative holds.  If $M_{\infty}$ is
noncompact the third alternative holds, and the theorem follows.
\end{proof}

\bibliographystyle{hamsplain}

\end{document}